\documentclass{cta-author}
\usepackage{algorithm}
\usepackage{algcompatible}
  \usepackage{natbib}
\usepackage{etex}
\usepackage{circuitikz}
\usetikzlibrary{shapes,arrows}
\usetikzlibrary{decorations.pathreplacing}
\usepackage{pgfplots}

\usepackage{graphics} 
\usepackage{epsfig} 
\usepackage{amsmath} 
\usepackage{mathdots}
\usepackage{amssymb}  
\usepackage{float}
\usepackage{graphicx}
\usepackage{placeins}
\usepackage{wrapfig}
\usepackage{subfig}

\usepackage{mathtools}

\usepackage{epstopdf}
\usepackage{fancyhdr}
\usepackage{booktabs,array}
\usepackage[output-decimal-marker={,}]{siunitx}
\usepackage{color}
\usepackage{empheq}

\newcommand{\al}[1]{\begin{align} #1 \end{align}}
\newcommand{\nn}{\nonumber}
\newcommand{\tr}{\mathrm{tr}}
\newcommand{\mz}{\color{black}}

\newtheorem{theorem}{Theorem}{}
{}
\newtheorem{remark}{Remark}{}
\newtheorem{proposition}{Proposition}{}
\newtheorem{lemma}{Lemma}{}
\begin{document}

\supertitle{Coupling MPC with Robust Kalman Filter}

\title{On the coupling of Model Predictive Control\\ and Robust Kalman Filtering}

\author{\au{Alberto Zenere$^{1}$}, \au{Mattia Zorzi$^{2}$}}

\address{\add{1}{Division of Automatic Control, Dept. of Electrical Engineering,
Link\"oping University, SE-58183, Link\"oping, Sweden. Email: alberto.zenere@liu.se}
\add{2}{Dipartimento di Ingegneria dell'Informazione, University of Padova, via Gradenigo 6/B, 35131 Padova, Italy. Email: zorzimat@dei.unipd.it}
}

\begin{abstract}
Model Predictive Control (MPC) represents nowadays one of the main methods employed for process control in industry. Its strong suits comprise a simple algorithm based on a straightforward formulation and the flexibility to deal with constraints. On the other hand it can be questioned its robustness regarding model uncertainties and external noises. Thus, a lot of efforts have been spent in the past years into the search of methods to address these shortcomings. In this paper we propose a robust MPC controller which stems from the idea of adding robustness in the prediction phase of the algorithm while leaving the core of MPC untouched. 
More precisely, we consider a robust Kalman filter that has been recently introduced and we further extend its usability to feedback control systems.
Overall the proposed control algorithm allows to maintain all of the advantages of MPC with an additional improvement in performance and without any drawbacks in terms of computational complexity. To test the actual reliability of the algorithm we apply it to control a servomechanism system characterized by nonlinear dynamics.
\end{abstract}

\maketitle

 \section{Introduction}

  Model Predictive Control (MPC), also referred to as Receding Horizon Control, is widely adopted in industry and represents the technology of choice to deal with the majority of constrained control problems \cite{Mayne20142967,bib:surveyMPC,camacho2012model,iet1,iet6}. MPC stems from the idea of employing a model that explains the relations among the variables of the plant to be controlled, which is then used to predict the future outputs. At each sampling time, a cost function $J_t$ comprising the tracking error (namely the difference between the desired trajectory and the predicted output) and the actuator efforts is minimized.
	The result of this procedure is an optimal sequence of future control moves which is applied according to a receding horizon philosophy: at time $t$ only the first input of the optimal command sequence is actually applied to the plant. The remaining optimal inputs are discarded, and a new minimization control problem is solved at time $t + 1$.

    It is important to keep in mind that the employed model is always an approximation of the actual process and the more accurate the approximation is, the better the MPC performance is. On the other hand there exists a trade-off between model accuracy and complexity of the optimization: the simpler the model is (and thus less accurate), the easier is solving the optimization problem. Therefore the idea is to construct a prediction model that is very simple but at same time representative enough to capture the main dynamical relations. Accordingly, a fundamental question about MPC regards its robustness with respect to model uncertainty and noise.
    Several different approaches can be found in literature to take into account model uncertainties when determining the optimal control sequence, \cite{CALAFIORE2013,bib:RMPCH,bib:surveyMPCrob2,bib:surveyMPCrob,KOTHARE19961361,bib:RMPCimp,iet2,iet3,iet4}.

	In practice however, almost none of the methods presented is adopted and it is preferred instead an {\em ad hoc} MPC tuning \cite[Sec. 1]{MACIEJOWSKI200922}, \cite[Sec. 10]{bib:surveyMPCrob2}. It is crucial to keep in mind that the user is looking for a controller that could actually be implemented in practice. It is then important to provide a control system that is feasible in terms of computational complexity and it is also easy to implement, even without an in-depth knowledge of the mathematics behind it.

 This paper follows exactly this direction. We propose a robust MPC hinging on the fundamental observation that the predictions used in MPC heavily rely on the accuracy of the employed state-space model. Hence the idea is to consider the usual MPC equipped with a robust Kalman filter {\mz whose predictions take into account the fact that the employed model is just an approximation of an accurate but complex (known or unknown) model of the actual process.}  More specifically in this paper we explore the use of the robust Kalman filter proposed in \cite{LN2013}, see also \cite{bib:rk,LN2004,HANSEN_SARGENT_2007,ROBUSTNESS_HANSENSARGENT_2008}. According to this approach, all possible incremental models belong to a ball which is formed by placing a bound on the {\em Kullback-Leibler} divergence (or possibly other divergences, \cite{BETA,DUAL}) between the actual and the nominal incremental model. This bound represents the tolerance accounting the mismodeling budget. Then, the robust filter is obtained by minimizing the mean square error according to the least favorable model in this ball.

 It is important to note that this robust Kalman filter was developed in the context of systems without inputs, \cite[Section III]{LN2013}. 
 Accordingly, in this paper we {\mz consider this robust filter} to the more general case of controlled feedback systems. Remarkably, the filter still admits a Kalman-like structure leading to a simple implementation of the corresponding MPC algorithm and allowing a reasonably low computational burden. 
 To test the effectiveness of the robust MPC we apply it to control a servomechanism system characterized by nonlinear dynamics. In particular, we compare its performance with standard MPC (i.e. MPC equipped with the Kalman filter) and MPC equipped with the risk-sensitive Kalman filter, \cite{bib:risk,speyer1998,HASSIBI_SAYED_KAILATH_BOOK,YOON_2004}.

    The outline of this paper is the following.
	In Sec.~\ref{sec:formulation} we briefly review classic MPC formulation. In Sec.~\ref{MPCproposed} we present the main concepts behind the considered robust Kalman filter as well as its extension to systems with feedback and thus its applicability combined with MPC. In Sec.~\ref{sec:DC} we introduce the physics of the servomechanism system that we considered and we show different simulations to attest the efficiency of the control system, lastly in Sec.~\ref{sec:conclusions} the conclusions are drawn.

\section{MPC Formulation} \label{sec:formulation}
\subsection{Standard MPC}

We briefly review the usual MPC formulation \cite[Chapter 2]{bib:macie}. We consider the discrete-time state-space model
\begin{align}
\Sigma\,:\,\begin{cases}
x_{t+1}& = Ax_{t} + Bu_{t}+Gv_t \\
\quad	y_{t}& = Cx_{t}+Dv_t
\end{cases}
\label{eq:model}
\end{align}
where $x_t\in\mathbb{R}^n$, $u_t\in\mathbb{R}^q$, $y_t\in\mathbb{R}^p$, $v_t\in\mathbb{R}^m$, denote the state, the input, the output and unmeasured noise, respectively. {\mz We assume that $x_0$ is Gaussian distributed with mean $\bar x_0$ and covariance matrix $\bar P_0$ positive definite.}
For simplicity in what follows we assume that $GD^T=0$, that is the state noise and the observation noise are independent.
Assume that our task is for the output $y_t$ to follow a certain reference signal $r_t$.
To this end, we introduce  the following quadratic cost function:
\begin{equation}
J_t (\textbf{u}_t,\Sigma)= \sum_{k=1}^{H_p}\Vert \hat{y}_{t+k|t}-r_{t+k}\Vert^2_{Q_k} + \sum_{k=0}^{H_u-1} \Vert \Delta \hat{u}_{t+k|t} \Vert^2_{R_k}.
\label{eq:costo}
\end{equation}
Here $\hat{y}_{t+k|t}$ represents the prediction of $y_{t+k}$ at time $t$ with $k>0$;  $\Delta\hat{u}_{t+k|t}:=\hat{u}_{t+k|t}-\hat{u}_{t+k-1|t}$ is the predicted variation of the input from time $t+k-1$ to $t+k$ and  $\textbf{u}_t := [u^T_{t|t}\ .\ .\ .\ u^T_{t+H_u-1|t}]^T$. 
The prediction and control horizons have length $H_p$ and $H_u$, respectively, and $H_u\leq H_p$. 
Lastly, $Q_k\in\mathbb{R}^{p\times p}$ and $R_k\in\mathbb{R}^{q\times q}$ indicate the weight matrices for the output prediction errors at time $t+k$ and for the predicted variations of the input at time $t+k$.
According to the receding horizon strategy, the control input $u_{t|t}$ to apply to $\Sigma$ at time $t$ is extracted
from $\textbf{u}_t$ which is solution to the following open-loop optimization problem
\al{ \label{optpb} u_{t|t}= \underset{\textbf{u}_t}{\mathrm{argmin}} \,J_t(\textbf{u}_t,\Sigma).}
The well known solution is, \cite[Chapter 3]{bib:macie}:

\al{ \label{utt}
u_{t|t} =
[I_q\ 0\ 0\ \cdots] (\Theta^T Q \Theta + R)^{-1}\Theta^T Q(\textbf{r}_t - \Psi \hat x_{t|t})
}
where
\begin{align*}
\Psi =&\left[\begin{array}{cccc} (CA)^T & (CA^2)^T & \ldots & (CA^{H_p})^T \end{array}\right]^T\\
\Theta =&
\begin{bmatrix}
CB & 0_{n\times q} & \cdots & 0_{n\times q}\\
CAB & CB & \cdots & 0_{n\times q}\\
\vdots && \ddots & \vdots  \\
CA^{H_p-1}B && \cdots  & CA^{H_p-H_u-1}B
\end{bmatrix}\\	
 Q=&\operatorname{diag}(Q_1,\cdots,Q_{H_p})\nn\\
  R=&\operatorname{diag}(R_0,\cdots,R_{H_u-1})\nn\\
\textbf{r}_t =&\left[\begin{array}{cccc} r_{t+1}^T & r_{t+2}^T & \ldots & r_{t+H_p}^T \end{array}\right]^T
\end{align*}
and $\hat x_{t|t}$ denotes the estimate of $x_t$ at time $t$. The latter is typically computed by using the Kalman filter
\al{ \label{KF1}\hat x_{t|t}&= \hat x_{t|t-1} +L_t (y_{t}-C \hat x_{t|t-1})\\
\label{KF2}\hat{x}_{t+1|t} &=  A\hat{x}_{t|t-1} + K_t(y_{t}-C\hat{x}_{t|t-1})+Bu_{t}  \\
 L_t &= {\mz P_tC^T(CP_tC^T + DD^T)^{-1},\;\; K_t=AL_t}\\
\label{KF4}P_{t+1} &= AP_tA^T - K_t(C P_tC^T+DD^T)K_t^T+GG^T  }
{\mz where the initial conditions are $\hat x_{0|-1}=\bar x_0$ and $P_0=\bar P_0$.}
  The resulting MPC law is outlined by Algorithm \ref{algoMPC}.
 \begin{algorithm}
\caption{MPC Law}
\label{algoMPC}
\begin{algorithmic}[1] \small
\STATE Collect the new data $y_t$
\STATE $L_t= P_t C^T(CP_t C^T+DD^T)^{-1}$
\STATE $\hat x_{t|t}=\hat x_{t|t-1} +L_t (y_{t}-C \hat x_{t|t-1})$
\STATE $u_{t|t} = [I_q\ 0\ 0\ \cdots] (\Theta^T Q \Theta + R)^{-1}\Theta^T Q(\textbf{r}_t - \Psi \hat x_{t|t})$
\STATE  Apply $u_{t|t}$ to the system
\STATE $K_t=A L_t $
\STATE $P_{t+1}= AP_tA^T-K_t(CV_tC^T+DD^T)K_t^T+GG^T$
\STATE $\hat x_{t+1|t}=A\hat{x}_{t|t-1} + K_t(y_{t}-C\hat{x}_{t|t-1})+Gu_{t} $
\STATE $t\leftarrow t+1$
\STATE Go to 1
\end{algorithmic}
\end{algorithm}
 Note that the estimator $\hat x_{t|t}$ is computed by assuming to know the actual underlying model $\Sigma$. In Section \ref{subsec_rob} we address the situation in which the actual model is different from $\Sigma$.

 It is worth noting that the optimization problem (\ref{optpb}) is usually considered with constraints, such as model uncertainty constraints and stability constraints. Such constraints are relatively easy to embed in (\ref{optpb}), on the other hand the price to pay is that the corresponding problem does not admit a closed form solution. Accordingly, the computational burden of Step 4 in Algorithm \ref{algoMPC} increases. In what follows we shall continue to consider the unconstrained MPC 
 because, as we will see, we will embed the model uncertainty in a different way.

 \subsection{Robust MPC} \label{subsec_rob}
 Assume that the actual model, denoted by $\tilde \Sigma$, is unknown and different from the nominal one, denoted by $\Sigma$.
It is then reasonable to assume that we are able to describe this uncertainty, that is we can characterize a set of models $\mathcal{S}$ for which $\tilde \Sigma\in\mathcal{S}$. In the robust MPC formulation the optimization problem (\ref{optpb}) is usually substituted by the mini-max problem \cite[Section 6]{bib:surveyMPCrob2}:
\al{ \label{minimaxpb}u_{t|t}= \underset{\textbf{u}_t}{\mathrm{argmin}} \,\underset{\tilde \Sigma\in\mathcal{S}}{\mathrm{max}} \,J_t(\textbf{u}_t,\tilde \Sigma).}
The latter is sometimes rewritten as a constrained MPC problem.
It is worth noting that solving a mini-max problem (or a constrained MPC problem) is computationally more demanding than solving a min problem. Many different uncertainty descriptions have been proposed in the literature such as impulse response uncertainty \cite{bib:RMPCimp}, structured feedback uncertainty, \cite{KOTHARE19961361}, polytopic uncertainty, \cite{iet5}, disturbances uncertainty, \cite{bib:RMPCH}, probabilistic uncertainty, \cite{LYGEROS_2011}, and Gaussian processes for modeling the underlying system, \cite{iet3}. Finally, in \cite{yang2015risk} it has been proposed a robust MPC wherein the cost function is an exponential-quadratic cost over the state distribution. In this way, large errors are severely penalized.

  \section{Robust MPC Proposed}\label{MPCproposed}
  As we already noticed, the standard MPC relies on the the assumption that the actual underlying model is known and thus the Kalman filter (\ref{KF1})-(\ref{KF4}) is designed on it. This assumption, therefore, could deteriorate the performance of MPC when the actual model is different from the nominal one. We propose
  a robust MPC which stems from the idea of building a control system consisting of MPC on one hand, but equipped with a robust state estimator, that takes into account possible differences between the actual and the nominal model, on the other.
In contrast with the usual robust MPC formulation, which is typically based on the mini-max Problem (\ref{minimaxpb}), we consider two independent optimization problems:
\begin{itemize}
\item Robust estimation problem: we want to find a robust estimate $\hat x_{t|t}$ of $x_t$, independently of the the fact that it will be next used to determine the optimal control input $u_{t|t}$. As we will see in Section \ref{sec:Rkalman} this problem is a mini-max problem itself, but its solution gives a robust filter obeying a Kalman-like recursion;
\item Open-loop optimization problem: assuming to have $\hat x_{t|t}$, we want to determine the optimal control input $u_{t|t}$, i.e. it coincides with Problem (\ref{optpb}).
\end{itemize}

  \subsection{Robust Kalman filter with input} \label{sec:Rkalman}
Recently, a new robust Kalman filter has been proposed \cite{LN2013,bib:rk,LN2004,bib:RKweiner} which represents the generalization of the risk-sensitive filters \cite{speyer1998,bib:risk}.

Assume that the nominal model $\Sigma$ is in the form (\ref{eq:model}), where $v_{t}\in\mathbb{R}^m$ is white Gaussian noise (WGN) with E[$v_{t}v^T_{s}$] = $\delta_{t-s}I_m$ ($\delta_t$ represents the Kronecker delta function). Moreover, the noise $v_t$ is independent of the initial state $x_0$, whose nominal distribution is given by $
f_0(x_0)\sim\mathcal{N}{\mz(\bar{x}_0,\bar P_0)}$. We introduce the random vector $z_t=[\,x_{t+1}^T\,\ y_t^T\,]^T$. At time $t$ the model $\Sigma$ is completely described by the conditional
 probability density of $z_t$ given the measurements $Y_{t-1}:=[\,y_0^T\,y_1^T\,  \ldots \,y_{t-1}^T\,]^T$, denoted by $
{f}_t(z_t|Y_{t-1})$. Note that by construction $ f_t(z_t|Y_{t-1})$  is Gaussian. Let $ \tilde f_t(z_t|Y_{t-1})$ be the conditional probability density of the actual underlying model $\tilde \Sigma$ at time $t$ and assume that it is Gaussian. The discrepancy between $\tilde \Sigma$ and $\Sigma$ at time $t$ can be measured through the Kullback-Leibler divergence
\al{\label{def_D_KL}\mathcal{D}(\tilde f_t, f_t)=\int_{\mathbb{R}^{n+p}} \tilde f_t(z_t|Y_{t-1}) \log\left(\frac{\tilde f_t(z_t|Y_{t-1})}{f_t(z_t|Y_{t-1})}\right) \mathrm{d}z_t.}
Thus, we define as the set of all allowable models at time $t$:
\begin{equation}
\mathcal{S}_{t} := \{ \tilde{f}_t(z_t|Y_{t-1})\ |\ \mathcal{D}(\tilde{f}_t , {f}_t) \le c \}
\label{eq:palla}
\end{equation}
that is, at time $t$ the actual model $\tilde \Sigma$ belongs to a ball  about $\Sigma$. The radius of this ball, denoted by $c$, represents the allowable model tolerance and must be fixed {\em a priori}.  Let $\mathcal{G}_t$ denote the class of estimators $g_t$ with finite second-order moments with respect to any probability density $\tilde{f}_{t}(z_{t}|Y_{t-1}) \in \mathcal{S}_{t}$.

We define as robust estimator of $x_{t+1}$ given $Y_{t}$ the solution to the following mini-max problem \begin{equation}
\hat x_{t+1|t} = \underset{g_{t}\in \mathcal{G}_{t}}{\operatorname{argmin}}\ \underset{\tilde{f}_{t}\in {\mathcal{S}}_{t}}{\operatorname{max}}\ \mathbb{E}_{\tilde f_t}[\| x_{t+1}-g_t(y_t)\|^2|Y_{t-1}]
\label{eq:minimax}
\end{equation}
where
\al{
\mathbb{E}_{\tilde f_t} & [\| x_{t+1}-g_t(y_t)\|^2|Y_{t-1}]\nn\\ & = \int_{\mathbb{R}^{n+p}}\ \Vert x_{t+1} - g_t(y_t) \Vert^2 \tilde{f}_t(z_t|Y_{t-1})dz_t\nn
} is the mean square error of the estimator with respect to the actual model $\tilde \Sigma$. Finally, the robust estimator of $x_t$ given $Y_{t}$, i.e. $\hat x_{t|t}$, is defined as the minimum mean square error based on (\ref{eq:model}) propagating $\tilde f_t(z_t|Y_{t-1})$.

In \cite{LN2013}  the solution of the mini-max problem (\ref{eq:minimax}) is characterized in the context of systems without inputs.
However the generalisation of this result to the more general case of feedback control systems is not trivial and needs to be properly addressed.
In the following theorem we consider a system with an input composed of a deterministic part plus a linear combination of the past values of the output, see Equation (\ref{utt}), and we show that the robust filter (\ref{eq:minimax}) still obeys to a Kalman-like recursion.
\begin{theorem} \label{teo1}
Let $u_t=a_t(Y_{t})+b_t(r_t)$, where $a_t(Y_{t})$ is a linear function of $Y_{t}$ and $b_t(r_t)$ a function of  a deterministic process $r_t$.
Then, the robust filter solution to (\ref{eq:minimax}) obeys to the recursion (\ref{KF1})-(\ref{KF2}) where
\al{\label{RK1bisa}V_t& =(P_t^{-1}-\theta_t I)^{-1}\\
\label{RK1}L_t&= V_t C^T(CV_t C^T+DD^T)^{-1},\;\; K_t=A L_t  \\
\label{RK1bis}P_{t+1}&= AV_tA^T-K_t(CV_tC^T+DD^T)K_t^T+GG^T ;}
$\theta_t>0$ is the unique solution to
\al{\label{RK2}-\log\det(I-\theta_t P_t)^{-1}+\tr[(I-\theta_t P_t)^{-1}-I]=c,} {\mz and the initial conditions are $\hat x_{0|-1}=\bar x_0$,  $P_0=\bar P_0$.}
Let $e_t=x_t-\hat x_{t|t-1}$ denote the prediction error of $x_t$ under the least-favorable model (solution of (\ref{eq:minimax})), then $e_t$ is Gaussian with zero mean and covariance matrix $V_t$.
\end{theorem}

\begin{remark} The right-hand side of (\ref{RK2}) is strictly convex and monotone increasing in $\theta$, \cite{bib:RSconvergence}. As a consequence, the solution can be easily found by using the bisection method, see Algorithm \ref{algotheta}.  
\end{remark}

\begin{algorithm}
\caption{Bisection algorithm for computing $\theta_t$}
\label{algotheta}
\begin{algorithmic}[1] \small
\STATE \textbf{Input:} $\theta_1=\varepsilon$, $\theta_2=\lambda^{-1}-\varepsilon$ where $\lambda$ max. eigenvalue of $P_t$, $\varepsilon>0$ sufficiently small
\STATE  \textbf{Repeat}
\STATE \quad $\theta_{new}=\frac{\theta_1+\theta_2}{2}$
\STATE  \quad \textbf{If} the right hand side of (\ref{RK2}) evaluated for $\theta=\theta_{new}$
is smaller than $c$
\STATE \quad \textbf{then}
\STATE \quad \quad $\theta_1=\theta_{new}$
\STATE \quad \textbf{else}
\STATE \quad \quad $\theta_2=\theta_{new}$
\STATE \quad \textbf{endif}
\STATE \textbf{Until} $|\theta_1-\theta_2|\leq \varepsilon$
\end{algorithmic}
\end{algorithm}

\begin{proposition} \label{stabR} Consider the robust filter (\ref{eq:minimax}) where $(A,G)$ is reachable and $(A,C)$ is observable. For $c>0$ in (\ref{eq:palla}) sufficiently small then $P_t \rightarrow \bar P$,
$V_t\rightarrow  \bar V$, $\theta_t \rightarrow  \bar \theta$, $L_t \rightarrow  \bar L$, $K_t\rightarrow  \bar K $ as $t\rightarrow  \infty$, that is the filter gain converges. Moreover, $\bar P>0$, $\bar V>0$, $\bar \theta>0$ and $\bar K$ is such that matrix $A- \bar KC$ is Schur stable meaning that in the steady state the prediction error $e_t$ is with bounded covariance matrix.     
\end{proposition}

{\mz  The choice of the tolerance parameter $c$ is a delicate step. In the design phase, the user typically has: (i) an accurate and complex model of the actual process; in practice this model can be understood  as an exact description of the actual process so that $\tilde\Sigma$ is known and coincides with this model; (ii) a simple model $\Sigma$ of the actual process of the form (\ref{eq:model}). Since $\tilde\Sigma$ and $\Sigma$ are known, then we can compute $\mathcal D(\tilde f_t,f_t)$ defined in (\ref{def_D_KL}).  Therefore, as an initial guess the parameter $c$ can be chosen as the average of $\mathcal D(\tilde f_t,f_t)$ over the interval of interest $t=1 \ldots t_{MAX}$.}

The constraint $\tilde f_t\in\mathcal{S}_t$, i.e. $\mathcal{D}(\tilde f_t, f_t)\leq c$, in (\ref{eq:minimax}) can be made soft by adding a penalty term in the objective function
 \al{
&\hat x_{t+1|t}=\nn\\  & \; \;\underset{g_t\in\mathcal{G}_t}{\operatorname{ argmin\ }}\underset{\tilde{f}_t\in\mathcal{V}_t}{\operatorname{ max\ }}  \mathbb{E}_{\tilde f_t}[\| x_{t+1}-g_t(y_t)\|^2|Y_{t-1}] - \bar \theta^{-1}\mathcal{D}(\tilde{f}_t , f_t)
\label{eq:rs2}
}
 where $\mathcal{V}_{t} := \{ \tilde{f}_t(z_t|Y_{t-1})\ |\ \mathcal{D}(\tilde{f}_t , {f}_t) \le \infty \}$ and  $\bar \theta>0$ is fixed {\em a priori}. It is worth noting that the solution to (\ref{eq:rs2}) is the risk-sensitive filter with risk sensitivity parameter $\bar \theta$, see
\cite{bib:risk}.

\begin{theorem}\label{teo2}
Let $u_t=a_t(Y_{t})+b_t(r_t)$, where $a_t(Y_{t})$ is a linear function of $Y_{t}$ and $b_t(r_t)$ a function of a deterministic process $r_t$.
Then, the risk-sensitive filter solution to (\ref{eq:rs2}) obeys to recursion (\ref{KF1})-(\ref{KF2}) where
\al{\label{eqVrs}V_t& =(P_t^{-1}-\bar \theta I)^{-1}\nn\\
L_t&= V_t C^T(CV_t C^T+DD^T)^{-1},\;\; K_t=A L_t  \nn\\
P_{t+1}&= AV_tA^T-K_t(CV_tC^T+DD^T)K_t^T+GG^T } 
{\mz and the initial conditions are $\hat x_{0|-1}=\bar x_0$,  $P_0=\bar P_0$ positive definite.}
Let $e_t=x_t-\hat x_{t|t-1}$ denote the prediction error of $x_t$ under the least-favorable model (solution of (\ref{eq:rs2})), then $e_t$ is Gaussian with zero mean and covariance matrix $V_t$.
\end{theorem}

\begin{remark}
The risk-sensitive filter of Theorem \ref{teo2} is not well defined for any $P_0$, on the other hand it is possible to impose some conditions on $P_0$
to guarantee the evolution of $P_t$ is on the cone of the positive definite matrices, \cite{levy2013contraction}. To avoid such a situation it is possible to substitute
(\ref{eqVrs}) with \al{\label{eqVrs2}V_t=F_t \exp(\bar\theta F_t^TF_t)F_t^T}
where $F_t$ is the Cholesky decomposition of $P_t$, i.e. $P_t=F_tF_t^T$. Relation (\ref{eqVrs2}) is given by solving the mini-max problem (\ref{eq:rs2})
 by replacing $\mathcal{D}(\tilde f_t, f_t)$ with the $\tau$-divergence, \cite{bib:RKdiv}, between $\tilde f_t$ and $f_t$ with parameter $\tau=1$, \cite{bib:rk}. \end{remark}

\begin{proposition} \label{stabRS} Consider the risk-sensitive filter (\ref{eq:rs2}) where $(A,G)$ is reachable and $(A,C)$ is observable.  For $	\bar \theta>0$ sufficiently small then $P_t \rightarrow \bar P$,
$V_t\rightarrow  \bar V$, $L_t \rightarrow  \bar L$, $K_t\rightarrow  \bar K $ as $t\rightarrow  \infty$. Accordingly, the filter gain converges. Moreover, $\bar P>0$, $\bar V>0$, and $\bar K$ is such that matrix $A- \bar KC$ is Schur stable meaning that in the steady state the prediction error $e_t$ is with bounded covariance matrix.     
\end{proposition}

{\mz To recap, in the presence of an input $u_t=a_t(Y_t)+b_t(r_t)$
we have modified the robust Kalman filter and the risk-sensitive 
 filter by adding the term $Bu_t$ in (\ref{KF2}). Even if the modification is very intuitive, it is not straightforward the fact that the modified filters are still robust (in the precise sense of \cite{bib:risk} and \cite{LN2013}, respectively). Theorem \ref{teo1} and Theorem \ref{teo2} showed that robustness is preserved in this case.}

\subsection{Robust MPC Law}
 The robust MPC law that we propose is constituted by two steps: first, we compute a robust estimate of $x_t$ given $Y_t$ according to (\ref{eq:minimax}); second, the
 optimal input control $u_{t|t}$ is computed as in (\ref{utt}). Since  $u_{t|t}$ is a function of $\{\mathbf{r}_t,\hat{x}_{t|t}\}$ and thus a function  of $\{\mathbf{r}_t, Y_t\}$, Theorem \ref{teo1} holds and the robust filter (\ref{eq:minimax}) obeys to a Kalman-like recursion.
 The resulting  robust MPC law is outlined in Algorithm \ref{algoMPCrob}.
\begin{algorithm}
\caption{Robust MPC Law}
\label{algoMPCrob}
\begin{algorithmic}[1] \small
\STATE Collect the new data $y_t$
\STATE Find $\theta_t$ s.t. $-\log\det(I-\theta_t P_t)^{-1}+\tr[(I-\theta_t P_t)^{-1}-I]=c$
\STATE $V_t =(P_t^{-1}-\theta_t I)^{-1}$
\STATE $L_t= V_t C^T(CV_t C^T+DD^T)^{-1}$
\STATE $\hat x_{t|t}=\hat x_{t|t-1} +L_t (y_{t}-C \hat x_{t|t-1})$
\STATE $u_{t|t} = [I_q\ 0\ 0\ \cdots] (\Theta^T Q \Theta + R)^{-1}\Theta^T Q(\textbf{r}_t - \Psi \hat x_{t|t})$
\STATE  Apply $u_{t|t}$ to the system
\STATE $K_t=A L_t $
\STATE $P_{t+1}= AV_tA^T-K_t(CV_tC^T+DD^T)K_t^T+GG^T$
\STATE $\hat x_{t+1|t}=A\hat{x}_{t|t-1} + K_t(y_{t}-C\hat{x}_{t|t-1})+Gu_{t} $
\STATE $t\leftarrow t+1$
\STATE Go to 1
\end{algorithmic}
\end{algorithm} To avoid the computation of 
$\theta_t$ in Step 2, we can approximate this robust filter with the risk-sensitive filter in Theorem \ref{teo2}. More precisely, since we know that $\theta_t\rightarrow \bar \theta$ as $t \rightarrow \infty$, we can approximate it with the risk-sensitive filter (\ref{eq:rs2}) with risk sensitivity parameter equal to $\bar \theta$.

\section{MPC of a Servomechanism System}
\label{sec:DC}

\begin{figure*}[t]
\resizebox {2\columnwidth} {!} {
\begin{circuitikz}
	
	\draw (0,3) to[V, v_=$V$] (0,0);
	\draw (0,3) to[R, i<^=$I_m$, l=$R$] (3,3);
	\draw (3,3) to[L, l=$L$] (4,3);
	
	\draw (4,3) -- (5,3);
	\draw (5,3) to[V, v_=$E_m$] (5,0);
	\draw (0,0) -- (5,0);
	
	\draw[fill=white] (4.85,0.85) rectangle (5.15,2.15);
	\draw[fill=white] (5,1.5) ellipse (.45 and .45);

	\draw[fill=white] (7,1.5)
	ellipse (.15 and 0.4);
	\draw (7.8,1.5) ellipse (.15 and 0.4);
	\draw[fill=white, color=white] (7.3, 1.1)
	rectangle (7.8, 1.9);
	\draw (7,1.1) -- (7.83,1.1);
	\draw (7,1.9) -- (7.83,1.9);	
	\draw (7.5,1.5) node {$J_m$};
	\draw (7.5,0.7) node {$\beta_m$};
	
	\draw[fill=black] (5.45,1.45) rectangle (7,1.55);
	
	\draw[line width=0.7pt,<-] (5.7,1) arc (-30:30:1);	
	\draw (6.3,1.8) node {$\theta_m$};	
	\draw (6.3,2.3) node {$T_m$};

	\draw[fill=black] (7.95,1.45) rectangle (10.2,1.55);
	
	\draw[line width=0.7pt,->] (8.1,1) arc (-30:30:1);	
	\draw[line width=0.7pt,->] (9.3,1) arc (-30:30:1);	
	\draw (8.6,1.8) node {$T_{f_m}$};	
	\draw (9.7,1.8) node {$T_s$};

	
%
%
%



	\draw[fill=black!50] (10.4,1.49)
	ellipse (.08 and 0.4);
	\draw[fill=black!50, color=black!50] (10.4,1.89)
	rectangle (10.2,1.09);
	\draw[fill=white] (10.2,1.49)
	ellipse (.08 and 0.4);
	\draw (10.2,1.89) -- (10.4,1.89);
	\draw (10.2,1.09) -- (10.4,1.09);

	\draw[fill=black!50] (10.4,0.40)
	ellipse (.13 and 0.67);
	\draw[fill=black!50, color=black!50] (10.4,1.07)
	rectangle (10.2,-0.27);
	\draw[fill=white] (10.2,0.40)
	ellipse (.13 and 0.67);
	\draw (10.2,1.07) -- (10.4,1.07);
	\draw (10.2,-0.27) -- (10.4,-0.27);


	\draw (9.6,0.4) node {$\rho$};
	
	\draw[fill=white] (11.5,0.4)
	ellipse (.15 and 0.4);
	\draw (12.3,0.4) ellipse (.15 and 0.4);
	\draw[fill=white, color=white] (11.8, 0)
	rectangle (12.3, 0.8);
	\draw (11.5,0) -- (12.33,0);
	\draw (11.5,0.8) -- (12.33,0.8);	
	\draw (12,0.4) node {$J_\ell$};
	\draw (12,-0.4) node {$\beta_\ell$};
	
	\draw[fill=black] (10.5,0.35) rectangle (11.5,0.45);
	
	\draw[fill=black] (12.45,0.35) rectangle (13.7,0.45);

	\draw[line width=0.7pt,<-] (12.6,0) arc (-30:30:1);	
	\draw[line width=0.7pt,->] (13.4,0) arc (-30:30:1);	
	\draw (13,0.7) node {$\theta_\ell$};
	\draw (13,1.2) node {$T_d$};		
	\draw (14,0.7) node {$T_{f_\ell}$};
	\end{circuitikz}
	}
\caption{Servomechanism System.}\label{fig:servomotore}
\end{figure*}
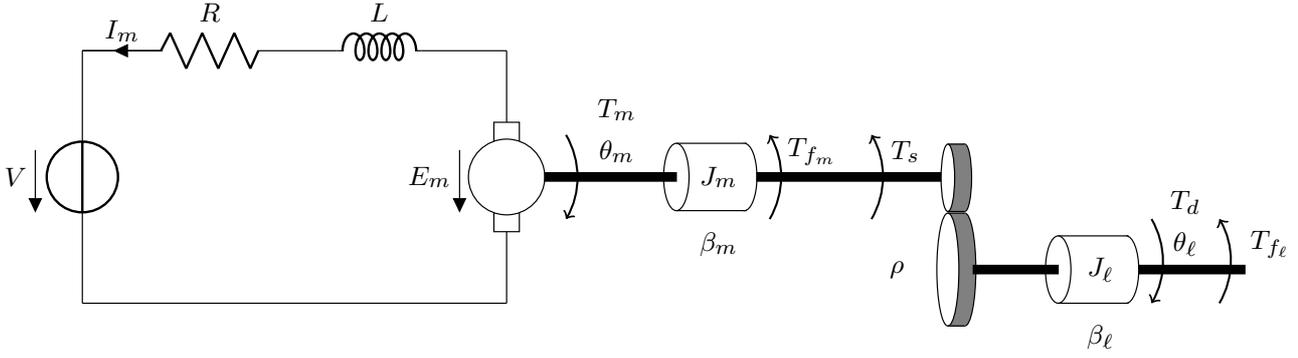

  To test the effectiveness of the robust MPC proposed in Section \ref{MPCproposed} we consider the servomechanism system presented in \cite{bib:DC}. It consists of a DC-motor, a gear-box, an elastic shaft and a load as depicted in Fig.\ref{fig:servomotore}. The input is the voltage applied $V$ and the output is the load angle $\theta_\ell$.
  Similar mechanisms are often found in the industry for a wide variety of applications.

  Since standard MPC is based on the assumption that the underlying model is linear, the general approach
  is to approximate the system with a linear model, possibly
neglecting the nonlinear dynamics.This approximation can be accurate enough as far as conventional control problems are concerned, however it may lead to intolerable errors when the DC motor operates at low speeds and rotates in two directions or when is needed a high precision control. In this situation, indeed, are significant the effects of the Coulomb and the deadzone friction which exhibit nonlinear behaviors in certain regions of operation. Next, we describe the simulation setup and we show the simulation results obtained by applying the standard MPC and the proposed robust MPC.

\subsection{Underlying Nonlinear Model}\label{sec:actual}
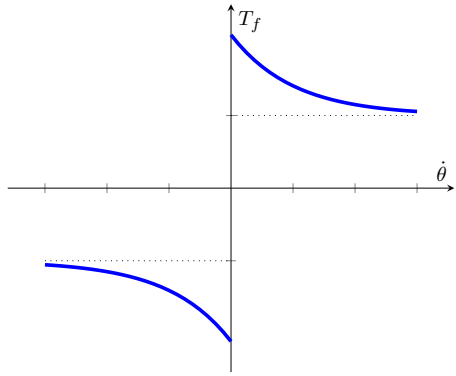
\begin{figure}
\begin{center}
\resizebox {0.7\columnwidth} {!} {
\begin{tikzpicture}
\begin{axis} [axis lines=middle,
enlargelimits,
xtick={},
ytick={-9,9},
xticklabels={},
yticklabels={},
xlabel=$\dot{\theta}$,ylabel=$T_f$]
\addplot[domain=0:15, blue, ultra thick] {9+10*exp(-0.2*x)};
\addplot[domain=-15:0, blue, ultra thick] {-9-10*exp(0.2*x)};
\addplot [dotted, domain=0:15] {9};
\addplot [dotted, domain=-15:0] {-9};
\end{axis}
\end{tikzpicture}
}
\end{center} \caption{Profile of the nonlinear friction model.} \label{fig_prof_nl_attr}
\end{figure}

  We now present the complete model including the nonlinear dynamics, that will represent the actual model in our simulations. The equations describing the physics of the system are:
	\begin{align*}
	J_\ell \ddot{\theta}_l &= \rho T_s - \beta_\ell\dot{\theta}_\ell - T_{f_\ell}(\dot{\theta}_\ell) \\
	J_m \ddot{\theta}_m &= T_m - T_s - \beta_m \dot{\theta}_m - T_{f_m}(\dot{\theta}_m) \\
	T_m &= K_t I_m \\
	V &= R I_m + L\dot{I}_m + E_m \\
	E_m &= K_t \dot{\theta}_m \\
	T_s &= \dfrac{k_\theta}{\rho}\left(\dfrac{\theta_m}{\rho}-\theta_\ell\right)
	\end{align*}
where
		$\theta_\ell$ denotes the load angle;
		$\theta_m$ the motor angle;
		$T_m$ the torque generated by the motor;
		$E_m$ the back electromotive force;
		$V$ the motor armature voltage;
		$I_m$ the armature current;
		 $T_s$ the torsional torque;
		$J_\ell$ the load inertia.
Moreover $T_{f_\ell}(\dot{\theta}_\ell)$ and  $T_{f_m}(\dot{\theta}_m)$ represent the Coulomb and the deadzone frictions on the load and on the motor, as described in detail in \cite{bib:attrito}:
\begin{align}
	T_{f_\ell}(\dot{\theta}_\ell) &= \alpha_{\ell 0}\operatorname{sgn}(\dot{\theta}_\ell) + \alpha_{\ell 1}e^{-\alpha_{\ell 2}\vert \dot{\theta}_\ell\vert}\operatorname{sgn}(\dot{\theta}_\ell)
	\label{eq:attrito_l}\\
	T_{f_m}(\dot{\theta}_m) &= \alpha_{m0}\operatorname{sgn}(\dot{\theta}_m) + \alpha_{m1}e^{-\alpha_{m2}\vert \dot{\theta}_m\vert}\operatorname{sgn}(\dot{\theta}_m)
	\label{eq:attrito_m}	
\end{align}
where the function $\operatorname{sgn}$ is defined as:
\[
	\operatorname{sgn}(x) =
	\begin{cases}
		1,\quad x > 0 \\
		0,\quad x = 0 \\
		-1,\quad x < 0.
	\end{cases}
\] The profile of the nonlinear friction model is depicted in Fig. \ref{fig_prof_nl_attr}

  The nominal values of the parameters of the servomechanism system are reported in Tab. \ref{tab:nom_values}.
  In practice, these values are not accurate because they are difficult to estimate. Accordingly,
 we introduce two kinds of possible parameter perturbations, $\varepsilon_{min}$ and $\varepsilon_{max}$,  expressing the percentage of the relative error that each nominal value can be affected by.
In the actual model every nominal parameter is perturbed in this way: if the nominal value is considered reliable enough it is perturbed by $\pm\varepsilon_{min}=\pm 5\%$, otherwise by $\pm\varepsilon_{max}=\pm 10\%$, see Tab.\ref{tab:real_values}.
\begin{table}[!h]
\processtable{Nominal parameters of the servomechanism system.\label{tab:nom_values}} 
{\begin{tabular}{ccl}
		\toprule
		\multicolumn{3}{c}{\textbf{Nominal parameters values}} \\
		\toprule
		\textit{Symbol} & \textit{Value(MKS)} & \textit{Meaning} \\
		\midrule	
		$L$         &  0       & Armature coil inductance \\
		$J_m$       &    0.5         &   Motor inertia      \\ 		
		$\beta_m$   &    0.1         &   Motor viscous friction coefficient     \\ 		
		$R$         &    20          &   Coil resistance of armature      \\ 		
		$K_t$       &    10          &   Motor constant      \\ 		
		$\rho$      &    20          &   Gear ratio      \\ 		
		$k_\theta$  &    1280.2      &   Torsional rigidity      \\ 		
		$J_\ell$ &    25    &   Load inertia      \\ 	
		$\beta_\ell$   &    25          &   Load viscous friction coefficient      \\ 												\\						\end{tabular} \hfill
	}{}
\end{table}

\begin{table}[!h]
		\processtable{Real parameters of the servomechanism system.\label{tab:real_values}}
{\begin{tabular}{ccl}
		\toprule
		\multicolumn{3}{c}{\textbf{Real parameters values}} \\
		\toprule
		\textit{Symbol} & \textit{Value(MKS)} & \textit{Meaning} \\
		\midrule
		$L$ & 0.8 & Armature coil inductance \\	
		$J_m$       &    0.5  (1+$\varepsilon_{max}$)        &   Motor inertia      \\ 		
		$\beta_m$   &    0.1  (1+$\varepsilon_{max}$)        &   Motor viscous friction coefficient     \\ 		
		$R$         &    20   (1+$\varepsilon_{min}$)        &   Resistance of armature      \\ 		
		$K_t$       &    10   (1 + $\varepsilon_{max}$)       &   Motor constant      \\ 		
		$\rho$      &    20   (1+$\varepsilon_{min}$)        &   Gear ratio      \\ 		
		$k_\theta$  &    1280.2 (1+$\varepsilon_{min}$)      &   Torsional rigidity      \\ 		
		${J}_\ell$ &    25 (1-$\varepsilon_{max}$)      &   Nominal load inertia      \\ 	
		$\beta_\ell$   &    25 (1+$\varepsilon_{max}$)          &   Load viscous friction coefficient      \\ 	
		$[\alpha_{\ell_0}\, \alpha_{\ell_1}\, \alpha_{\ell_2}]$   & [0.5\, 10\, 0.5]             &  Load Nonlinear friction parameters        \\ 								$[\alpha_{m_0}\, \alpha_{m_1}\, \alpha_{m_2}]$   & [0.1\, 2\, 0.5]             &      Motor Nonlinear friction parameters           \\ 	 \\																		
	\end{tabular}
}{}
\end{table}


\subsection{Linear Model for MPC} \label{secNomimod}
We now derive the nominal model. To obtain a linearized model of the servomechanism system we eliminate the nonlinear dynamics (\ref{eq:attrito_l})-(\ref{eq:attrito_m}) and set $L=0$.
  The dynamic equations resulting from these simplifications are:
\begin{align*}
	\begin{cases}
		J_\ell\ddot{\theta}_\ell &= \rho T_s - \beta_\ell \dot{\theta}_\ell \\
		J_m\ddot{\theta}_m &= T_m - T_s - \beta_m\dot{\theta}_m
	\end{cases}
\end{align*}
wherein we consider the nominal parameters in Tab. \ref{tab:nom_values}.
Defining as state vector  $x_t = [\,\theta_\ell\ \ \dot{\theta}_\ell\ \ \theta_m\ \ \dot{\theta}_m\,]^T$, we obtain a continuous-time state-space linear model of type

\begin{equation*}
\begin{cases}
\dot{x}_t = \bar A x_t + \bar Bu_t + \bar G \dot{w}_t \\
y_t = \bar Cx_t  + \bar D  \dot{w}_t \\
\end{cases}
\end{equation*}
where $u_t$ is the motor armature voltage, $y_t$ the load angle, $w_t$ is the normalized Wiener process, and matrices $\bar G$, $\bar D$ are chosen heuristically in such a way to compensate the approximations made before. Lastly this model was discretized with sampling time $T=0.1$ s obtaining in this way a model of type (\ref{eq:model}).

\subsection{Results}
\label{sec:results}
 In this section we want to show the improvement in performance brought by the MPC controller equipped with the robust Kalman filter in Theorem \ref{teo1} (R-MPC). As terms of comparison we consider two other controllers: standard MPC (S-MPC) and MPC equipped with the risk-sensitive filter in Theorem \ref{teo2} (RS-MPC).

\begin{figure}[!h]
	\centering
	\includegraphics[width=0.5\textwidth]{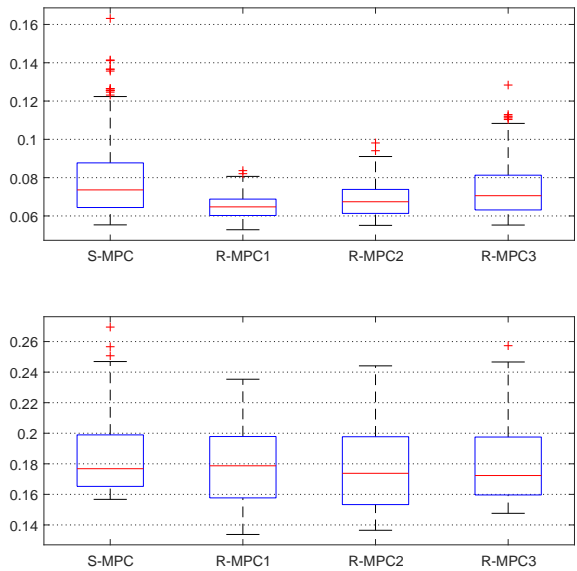}
	\caption{{\mz Boxplots of the MSE for S-MPC, R-MPC1, R-MPC2, R-MPC3 for different choices of the prediction and control horizons. The upper boxplot corresponds to $H_p = 10$ and $H_u = 8$ while the lower one to $H_p = 15$ and $H_u = 3$.} }
	\label{fig:box_horizons}
\end{figure}

More specifically we study their performances in response with the reference trajectory set to $r_t = \pi/2$ rad, $t\ge0$ with $y_0=0$ rad. The initial state $x_0$ is assumed to be zero mean and with variance equal to the variance of the state-process noise in (\ref{eq:model}). Regarding the parameters of MPC in the cost function (\ref{eq:costo}) we consider the weight matrices $Q_k = 0.1$, $R_k = 0.1$ and $H_p = 10$, $H_u = 3$.
{\mz In particular the effects of the choice of $H_p$ and $H_u$ are further analyzed in the last simulation, see Fig.\ref{fig:box_horizons}.}

\begin{figure}[!h]
	\centering
	\includegraphics[width=0.5\textwidth]{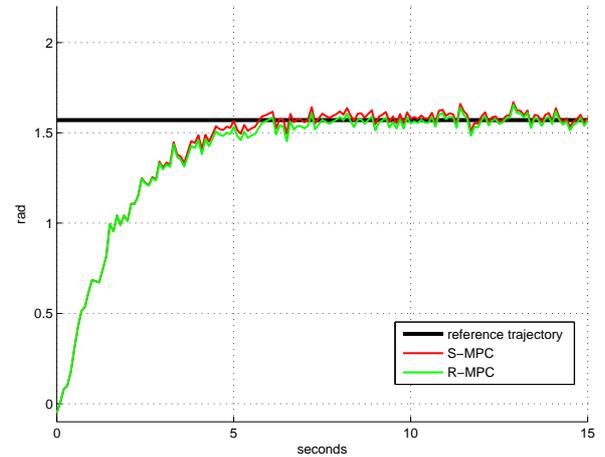}
	\caption{First simulation: load angle when the nominal and the actual model coincide.}
	\label{fig:1_output}
\end{figure}

\begin{figure}[!h]
	\centering
	\includegraphics[width=0.5\textwidth]{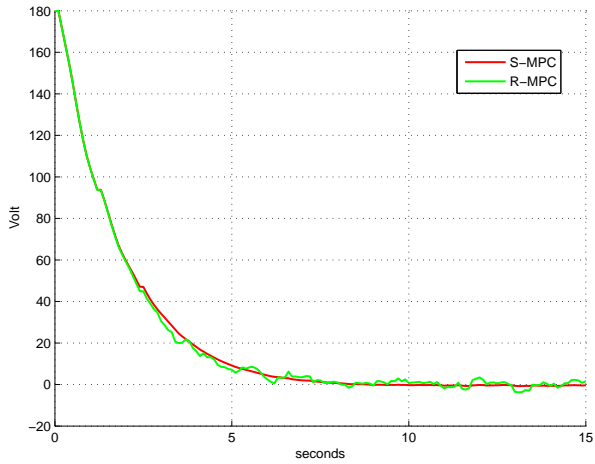}
	\caption{First simulation: voltage applied when the nominal and the actual model coincide.}
	\label{fig:1_input}
\end{figure}

In the first simulation we consider the ideal situation in which the actual model coincides with the one of Section \ref{secNomimod}, that is the actual and the nominal models coincide.
As can be seen in Fig.\ref{fig:1_output},  S-MPC can track the reference trajectory likewise R-MPC (with $c=10^{-1}$).
However, it is worth noting that the input applied by R-MPC is less smooth, see Fig. \ref{fig:1_input}. 
This is not a very surprising fact given that the robust filter is constructed under the idea of considering more uncertainties in the modelling. The performance of RS-MPC has been omitted because similar to the one of R-MPC.

\begin{figure}[!h]
	\centering
	\includegraphics[width=0.5\textwidth]{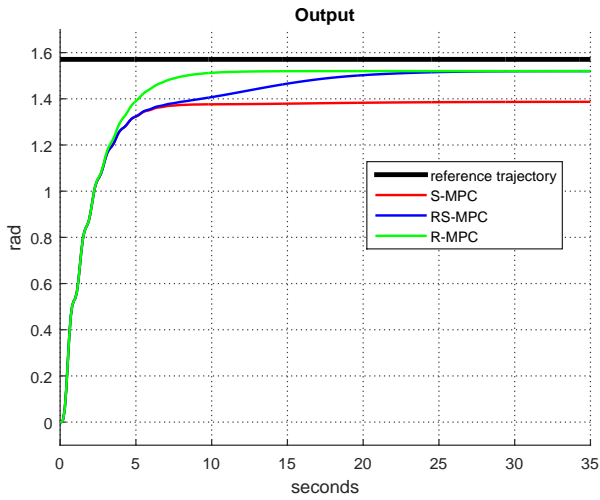}
	\caption{Second simulation: load angle when the actual model differs from the nominal one.}
	\label{fig:2_output}
\end{figure}

In the second simulation, the actual model is the one of Section \ref{sec:actual}. Fig. \ref{fig:2_output}
shows the load angle for S-MPC, R-MPC (with $c=10^{-1}$) and RS-MPC (with $	\bar \theta$ set to the value at steady state computed through R-MPC). 

\begin{figure}[!h]
	\centering
	\includegraphics[width=0.5\textwidth]{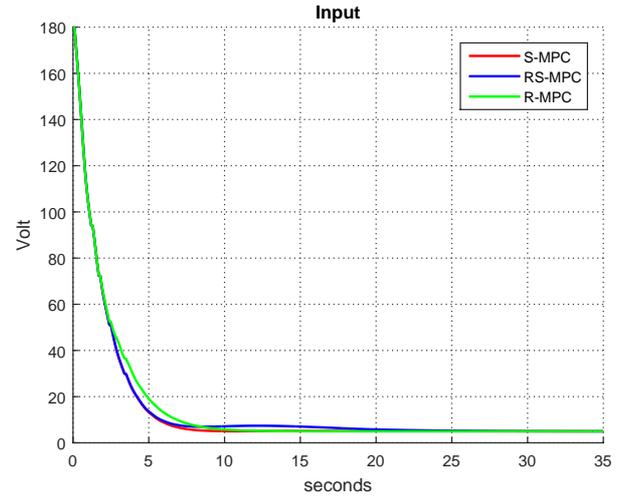}
	\caption{Second simulation: voltage applied when the actual model differs from the nominal one.}
	\label{fig:2_input}
\end{figure}

The differences between the three controllers are clearly visible.
In particular it is evident how R-MPC is able to provide an adequate control whereas RS-MPC and S-MPC do not. 
To give a quantitative analysis, in the first simulation we can see how after 5 seconds the error on the output is comprised within 5\% of the desired trajectory for all the controls (i.e. S-MPC, RS-MPC and R-MPC).
On the other hand, in the second simulation the error on the output is comprised within 5\% of the desired trajectory after 8 seconds for R-MPC, after 19 seconds for RS-MPC and the one of S-MPC never reaches this accuracy in 35 seconds.
Again, R-MPC is the unique control providing similar performances both in the first and in the second simulation.
On the other hand, R-MPC requires more energy (see Fig. \ref{fig:2_input} on the interval 3$\div$10 seconds) than S-MPC. Finally, RS-MPC constitutes an in-between solution in terms of performance.

\begin{figure}[!h]
	\centering
	\includegraphics[width=0.5\textwidth]{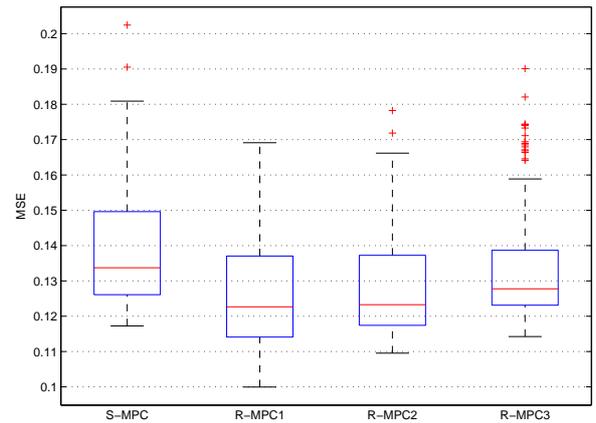}
	\caption{Boxplot of the MSE for S-MPC, R-MPC1, R-MPC2, R-MPC3: {\mz the central mark (in red) indicates the median, and the bottom and top edges of the blue box indicate the 25th and 75th percentiles.} }
	\label{fig:bloxplot}
\end{figure}

  The third simulation has the aim to evaluate the effect of the tolerance $c$ on the performance of R-MPC, {\mz see Fig. \ref{fig:bloxplot}}.
  More precisely we consider  R-MPC1 with $c=10^{-1}$, R-MPC2 with $c=10^{-2}$, and R-MPC3 with $c=10^{-3}$.  We consider a Monte Carlo experiment of 200 runs. At each run we generate the actual model where $\varepsilon_{min}$ and $\varepsilon_{max}$ in Tab. \ref{tab:real_values} are randomly chosen inside the interval, respectively,
[-10\%,10\%] and [-20\%,20\%] following a uniform distribution probability.
The unique exception is for the load inertia $J_\ell$ which is perturbed randomly with a relative error uniformly distributed inside the interval [-80\%,80\%]. In this way we consider the situation in which the servomechanism system is connected to different loads during its operational life. Then,  we apply the MPC controls and for each of them we compute the MSE of the output with respect to the reference signal during the first 20 seconds (in which the control with R-MPC usually approaches the steady state, as seen in Fig.\ref{fig:2_output}).

 As expected, R-MPC performs more accurately than S-MPC and, interestingly enough, 
  higher values of $c$ correspond to a diminished error of the MSE of R-MPC.
  On the other hand,   higher values of $c$  correspond to 
  high-energy control inputs.

{\mz As a last simulation, we analyze the influence of the length of the prediction and control horizons on the performance.
To this aim, we have carried out the third simulation for different choices of the horizons. 
As representative examples, in Fig. \ref{fig:box_horizons} we report the results corresponding to $H_p = 10, H_u = 8$ and $H_p = 15, H_u = 3$.
The first case shows that a large value for the control horizon significantly reduces the tracking error for all the MPC controllers. In particular, R-MPC1 provides the best performance while 
S-MPC provides the worst one. Conversely, a large value for the  prediction horizon (second case) is detrimental to the performance of the MPC controllers. In particular, all the MPC controllers perform in a similar way and R-MPC3 is slightly better than the others. Finally, we have noticed that the control action is very responsive with $H_u$ large and very slow with $H_p$ large.}


\section{Conclusions}
\label{sec:conclusions}

  Model Predictive Control is nowadays one of the primary methods to address process control and is indeed employed in a wide variety of applications. One main challenge that is now presented to the research community is to find a way to gather, in the same controller, the advantages of MPC with robustness properties.

  In this paper we explored the idea to build a controller which is constituted  by the usual MPC algorithm equipped with a robust Kalman filter.
  This appears as a very attractive option since it would meet the challenge above; in fact the design of MPC remains untouched with all of its benefits and the robust filtering can potentially compensate inaccuracies in modelling. Moreover, it would not deteriorate the computational cost from standard MPC.

  Thereafter we tested the controller taking into consideration a servomechanism system characterized by nonlinear dynamics. In order to obtain more realistic simulations, a margin of error on the values of the system parameters was introduced.
  To assess the capabilities of this controller, we evaluated its efficiency with respect to classic MPC and to MPC equipped with the risk-sensitive Kalman filter.

  Overall, the control system proved to be able to compensate errors in modelling significantly more effectively than the other two.
  In conclusion, the controller designed in this paper appears to be a viable solution to obtain all of the advantages of MPC while guaranteeing at the same time higher robustness.

\section{Appendix}
 \subsection{Proof of Theorem \ref{teo1}}
 To prove the statement we need the following result, see \cite{bib:RKweiner,LN2004}.
 \begin{lemma} \label{propo1}Let $z=[\,x^T\, y^T\,]$ be a random vector with nominal and actual probability density, respectively,
 $f\sim \mathcal{N}(m_z,K_z)$ and $\tilde f\sim \mathcal{N}(\tilde m_z,\tilde K_z)$. We conformably partition the mean and the covariance of $f$:
 \al{m_z=\left[\begin{array}{c} m_x\\ m_y\end{array}\right],\;\; K_z=\left[\begin{array}{cc} K_x & K_{x,y}\\ K_{y,x} & K_y\end{array}\right]} and likewise for $\tilde m_z$ and $\tilde K_z$.
Consider the mini-max problem
\al{\label{minimaxstatic}\underset{g\in \mathcal{G}}{\operatorname{min}}\ \underset{\tilde{f}\in\bar{\mathcal{S}}}{\operatorname{max}}\ \int_{\mathbb{R}^{n+p}} \| x-g(y)\|^2 \tilde f(z) \mathrm{d} z}
where $\mathcal{G}$ denotes the set of all estimators $g(y)$ with finite second order moments with respect to any $\tilde f\in\mathcal{S}$ and $\mathcal{S}=\{ \tilde f\,:\; \mathcal{D}(\tilde f, f)\leq c\}$. Then, the least favorable density $\tilde f$ solution to (\ref{minimaxstatic}) is such that
 \al{\tilde m_z=m_z,\;\; \tilde K_z=\left[\begin{array}{cc} \tilde K_x & K_{x,y}\\ K_{y,x} & K_y\end{array}\right].}
Let
\al{ P&=K_x-K_{x,y}K_y^{-1} K_{y,x} \nn\\
V&=\tilde K_x-K_{x,y}K_y^{-1} K_{y,x} }
be the nominal and the least favorable, respectively,  a posteriori error covariance matrix of $x$ given $y$. Then
$V=(P^{-1}-\theta I)^{-1}$
and $\theta>0$ is the unique solution to
\al{-\log\det(I-\theta P)^{-1}+\tr[(I-\theta P)^{-1}-I]=c.} The robust estimator solution to (\ref{minimaxstatic}) is the Bayes estimator
\al{g(y)=m_x+K_{x,y}K_y^{-1} (y-m_y).}
 \end{lemma}

In view of Lemma \ref{propo1}, it is worth noting that the estimation error $e:= x-g(y)$ under the least favorable model is Gaussian and with zero mean and covariance matrix $V$.

 We now proceed to prove the statement. Assume that the statement holds at time $t$, that is the robust estimator of
 $x_t$ given $Y_{t-1}$ is $\hat x_{t|t-1}=\mathbb{E}_{\tilde f_t}[x_t|Y_{t-1}]$ and the prediction error $e_t=x_t-\hat x_{t|t-1}$ is Gaussian with zero mean and covariance matrix $V_t$. We consider now the nominal model (\ref{eq:model}). Since $a_t(Y_t)$ is a linear function, we have \al{u_t&=\tilde a_t(Y_{t-1})+J_t y_t+b_t(r_t)\nn\\
&= \tilde a_t(Y_{t-1})+J_tC x_t+J_t D v_t+b_t(r_t)}
 where $\tilde a_t(Y_{t-1})$ is a linear function of $Y_{t-1}$ and $J_t$ a matrix.
Then, we have
\al{ f_t(z_t|Y_{t-1})\sim \mathcal{N}\left( \left[\begin{array}{c}\tilde A_t \\ C\end{array}\right] \hat x_{t|t-1}+
\left[\begin{array}{c} \tilde S_t  \\0\end{array}\right],K_{z_t}
\right)}
with
 \al{\label{def_K_zt}K_{z_t}&=\left[\begin{array}{cc} K_{x_{t+1}} & K_{x_{t+1},y_t}\\ K_{y_t,x_{t+1}} & K_{y_t}\end{array}\right]\nn\\
 & = \left[\begin{array}{c}\tilde A_t\\ C\end{array}\right] V_t \left[\begin{array}{cc}\tilde A_t^T  & C^T\end{array}\right]+
 \left[\begin{array}{c}\tilde G_t \\ D\end{array}\right]\left[\begin{array}{cc}\tilde G_t^T  & D^T\end{array}\right],}
 $\tilde S_t=B(\tilde a_t(Y_{t-1})+b_t(r_t))$, $\tilde A_t=A+BJ_tC $ and $\tilde G_t=G+BJ_tD$.
  Applying Lemma \ref{propo1}, we have that
 \al{ \tilde f_t(z_t|Y_{t-1})\sim \mathcal{N}\left( \left[\begin{array}{c}\tilde A_t \\ C\end{array}\right] \hat x_{t|t-1}+
\left[\begin{array}{c}\tilde S_t  \\0\end{array}\right],\tilde K_{z_t}
\right)}
where
\al{\tilde K_{z_t}&=\left[\begin{array}{cc} \tilde K_{x_{t+1}} & K_{x_{t+1},y_t}\\ K_{y_t,x_{t+1}} & K_{y_t}\end{array}\right].}
Accordingly, \al{\hat x_{t+1|t}&=\mathbb{E}_{\tilde f_t} [x_{t+1}|Y_{t-1}]\nn\\
\label{eq_filtro}&= \tilde A_t \hat x_{t|t-1}+\tilde S_t+K_{x_{t+1},y_t} K_{y_t}^{-1}(y_t-C \hat x_{t|t-1})} and $e_{t+1}=x_{t+1}-\hat x_{t+1|t}$ is Gaussian with zero mean and 
nominal and least favorable covariance matrix at time $t+1$, respectively, 
\al{ \label{Pt1}P_{t+1}&=K_{x_{t+1}}-K_{x_{t+1},y_t}K_{y_t}^{-1} K_{y_t,x_{t+1}}\\
V_{t+1}&=(P_{t+1}-\theta_{t}I)^{-1}\nn} where $\theta_t$ is the unique solution to (\ref{RK2}). By (\ref{eq_filtro}) and (\ref{def_K_zt}) we have
\al{\hat x_{t+1|t} &= \tilde A_t \hat x_{t|t-1}+\tilde S_t \nn \\
& \quad +(\tilde A_t V_t C^T+\tilde G_t D^T)K_{y_t}^{-1}(y_t-C \hat x_{t|t-1})\nn\\
&=\tilde A_t \hat x_{t|t-1}+\tilde S_t +(A V_t C^T+B J_t C V_t C^T+BJ_tDD^T)\nn \\
& \quad \times K_{y_t}^{-1}(y_t-C \hat x_{t|t-1})\nn\\
&=\tilde A_t \hat x_{t|t-1}+\tilde S_t +(A V_t C^T+B J_t K_{y_t})\nn \\
& \quad  \times K_{y_t}^{-1}(y_t-C \hat x_{t|t-1})\nn\\
&=A \hat x_{t|t-1}+B(\tilde a_t(Y_{t-1})+b_t(r_t) )+B J_t C \hat x_{t|t-1}\nn\\
& \quad
+A V_t C^TK_{y_t}^{-1}(y_t-C \hat x_{t|t-1})
+B J_t(y_t-C \hat x_{t|t-1})\nn\\
& = A \hat x_{t|t-1}+B(\tilde a_t(Y_{t-1})+b_t(r_t)+J_t y_t )\nn \\
& \quad +A V_t C^TK_{y_t}^{-1}(y_t-C \hat x_{t|t-1})\nn\\
&=A \hat x_{t|t-1}+Bu_t+A V_t C^TK_{y_t}^{-1}(y_t-C \hat x_{t|t-1})\nn
} which coincides with (\ref{KF2}) where $K_t$ has been defined in (\ref{RK1}). By (\ref{Pt1}) and (\ref{def_K_zt}) we have:
\al{P_{t+1}&= \tilde A_t V_t \tilde A_t^T+\tilde G_t \tilde G_t^T\nn \\
& \quad -(\tilde A_t V_t C^T +\tilde G_t D^T ) K_{y_t}^{-1} (\tilde A_t V_t C^T +\tilde G_t D^T )^T\nn\\
&= \tilde A_t V_t \tilde A_t^T+\tilde G_t \tilde G_t^T \nn \\
& \quad -
(A V_t C^T +B J_t C V_t C^T+B J_t D D^T ) \nn\\
& \quad \times  K_{y_t}^{-1} (A V_t C^T +B J_t C V_t C^T+B J_t D D^T )^T\nn\\
& =\tilde A_t V_t \tilde A_t^T+\tilde G_t \tilde G_t^T \nn \\
& \quad -
(A V_t C^T +B J_t K_{y_t})  K_{y_t}^{-1} (A V_t C^T +B J_t K_{y_t})^T \nn\\
&=  A V_t  A^T+B J_t C V_t A^T+A V_t C^T J_t^T B^T\nn \\
& \quad + B J_t C V_t C^T J_t^T B^T +GG^T +BJ_t DD^T J_t^T B^T\nn\\
& \quad -
A V_t C^T  K_{y_t}^{-1} A V_t C^T  -A V_t C^T   J_t^T B^T \nn \\
& \quad 
-B J_t C V_t A^T -B J_t K_{y_t} J_t^T B^T \nn\\
&=  A V_t  A^T  -
A V_t C^T  K_{y_t}^{-1} A V_t C^T +GG^T \nn} 
which coincides with (\ref{RK1bis}). Regarding the update equation (\ref{KF1}), it can be derived with the usual argumentations.
 \qed

   \subsection{Proof of Proposition \ref{stabR}}
It is sufficient to note that the recursion (\ref{RK1bisa})-(\ref{RK1bis}) is the one studied in \cite{convtau} with $\tau=0$. In \cite[Theorem 4.1]{convtau} it was shown that such a recursion converges provided that $c>0$ is sufficiently small. Moreover, matrix $A-\bar KC$ is Schur stable.\qed

 \subsection{Proof of Theorem \ref{teo2}}
  To prove the statement we need the following result, see \cite[Corollary 3.1]{bib:RKweiner} with $\tau=0$.
 \begin{lemma} Let $z=[\,x^T\, y^T\,]$ be a random vector with nominal and actual probability density, respectively,
 $f\sim \mathcal{N}(m_z,K_z)$ and $\tilde f\sim \mathcal{N}(\tilde m_z,\tilde K_z)$. We conformably partition the mean and the covariance of $f$:
 \al{m_z=\left[\begin{array}{c} m_x\\ m_y\end{array}\right],\;\; K_z=\left[\begin{array}{cc} K_x & K_{x,y}\\ K_{y,x} & K_y\end{array}\right]} and likewise for $\tilde m_z$ and $\tilde K_z$.
Consider the min problem
\al{\label{minimaxstatic}\underset{g\in \mathcal{G}}{\operatorname{min}}\ \int_{\mathbb{R}^{n+p}} \| x-g(y)\|^2 \tilde f(z) \mathrm{d} z-\theta^{-1} \mathcal D(\tilde f, f)}
where $\mathcal{G}$ denotes the set of all estimators $g(y)$ with finite second order moments with respect to any $\tilde f\in\mathcal{S}$ and $\mathcal{S}=\{ \tilde f\,:\; \mathcal{D}(\tilde f, f)\leq \infty\}$. Then, the least favorable density $\tilde f$ solution to (\ref{minimaxstatic}) is such that
 \al{\tilde m_z=m_z,\;\; \tilde K_z=\left[\begin{array}{cc} \tilde K_x & K_{x,y}\\ K_{y,x} & K_y\end{array}\right].}
Let
\al{ P&=K_x-K_{x,y}K_y^{-1} K_{y,x} \nn\\
V&=\tilde K_x-K_{x,y}K_y^{-1} K_{y,x} }
be the nominal and the least favorable, respectively,  a posteriori error covariance matrix of $x$ given $y$. Then
$V=(P^{-1}-\theta I)^{-1}$. The robust estimator solution to (\ref{minimaxstatic}) is the Bayes estimator
\al{g(y)=m_x+K_{x,y}K_y^{-1} (y-m_y).}
 \end{lemma}

 Using the above result, the proof of the Theorem is similar to the one of Theorem \ref{teo1}.\qed

  \subsection{Proof of Proposition \ref{stabRS}}
  It is sufficient to note that the recursion (\ref{eqVrs}) is the one studied in \cite{levy2013contraction}. In \cite[Theorem 5.3]{levy2013contraction} it was shown that such a recursion converges provided that $\bar \theta>0$ is sufficiently small. Moreover, matrix $A-\bar KC$ is Schur stable.\qed


\begin{thebibliography}{10}

\bibitem{iet6}
Peyman Bagheri, Peyman~Bagheri Ali, Ali~Khaki Sedigh, and Khaki Sedigh.
\newblock Analytical approach to tuning of model predictive control for
  first-order plus dead time models.
\newblock {\em IET Control Theory \& Applications}, 7(14):1806--1817, 2013.

\bibitem{speyer1998}
R.~N. Banavar and J.~L. Speyer.
\newblock Properties of risk-sensitive filters/estimators.
\newblock {\em IEEE Proceedings - Control Theory and Applications},
  145(1):106--112, Jan 1998.

\bibitem{bib:surveyMPCrob2}
A.~Bemporad and M.~Morari.
\newblock Robust model predictive control: a survey.
\newblock {\em {Lecture Notes in Control and Information Sciences}},
  245:207--226, 2007.

\bibitem{bib:DC}
A.~Bemporad and E.~Mosca.
\newblock Fulfilling hard constraints in uncertain linear systems by reference
  managing.
\newblock {\em Automatica}, 34(4):451--461, 1998.

\bibitem{bib:risk}
R.~K. Boel, M.~R. James, and I.~R. Petersen.
\newblock Robustness and risk-sensitive filtering.
\newblock {\em {IEEE Transactions on Automatic Control}}, 47(3):451--461, 2002.

\bibitem{CALAFIORE2013}
G.~C. Calafiore and L.~Fagiano.
\newblock Robust model predictive control via scenario optimization.
\newblock {\em IEEE Transactions on Automatic Control}, 58(1):219--224, Jan
  2013.

\bibitem{camacho2012model}
E.~F. Camacho and C.~Bordons.
\newblock {\em Model predictive control in the process industry}.
\newblock Springer Science \& Business Media, 2012.

\bibitem{iet3}
Gang Cao, Edmund M-K Lai, and Fakhrul Alam.
\newblock Gaussian process model predictive control of unknown non-linear
  systems.
\newblock {\em IET Control Theory \& Applications}, 11(5):703--713, 2017.

\bibitem{LYGEROS_2011}
D.~Chatterjee, P.~Hokayem, and J.~Lygeros.
\newblock Stochastic receding horizon control with bounded control inputs: A
  vector space approach.
\newblock {\em IEEE Transactions on Automatic Control}, 56(11):2704--2710, Nov
  2011.

\bibitem{HANSEN_SARGENT_2007}
L.~Hansen and T.~Sargent.
\newblock Recursive robust estimation and control without commitment.
\newblock {\em Journal of Economic Theory}, 136(1):1--27, 2007.

\bibitem{ROBUSTNESS_HANSENSARGENT_2008}
L.~Hansen and T.~Sargent.
\newblock {\em Robustness}.
\newblock Princeton University Press, Princeton, NJ, 2008.

\bibitem{HASSIBI_SAYED_KAILATH_BOOK}
B.~Hassibi, A.~Sayed, and T.~Kailath.
\newblock {\em Indefinite-Quadratic Estimation and Control-- A Unified Approach
  to $H^2$ and $H^{\infty}$ Theories}.
\newblock Soc. Indust. Appl. Math., Philadelphia, 1999.

\bibitem{bib:surveyMPCrob}
A.~A. Jalali and V.~Nadimi.
\newblock A survey on robust model predictive control from 1999-2006.
\newblock In {\em {Internation Conference on Computational Intelligence for
  Modelling Control and Automation and International Conference on Intelligent
  Agents, Web Technologies and Internet Commerce}}, 2006.

\bibitem{bib:attrito}
T.~Kara and I.~Eker.
\newblock Nonlinear modeling and identification of a \{DC\} motor for
  bidirectional operation with real time experiments.
\newblock {\em Energy Conversion and Management}, 45(7-8):1087 -- 1106, 2004.

\bibitem{KOTHARE19961361}
M.~V. Kothare, V.~Balakrishnan, and M.~Morari.
\newblock Robust constrained model predictive control using linear matrix
  inequalities.
\newblock {\em Automatica}, 32(10):1361 -- 1379, 1996.

\bibitem{LN2004}
B.~C. Levy and R.~Nikoukhah.
\newblock Robust least-squares estimation with a relative entropy constraint.
\newblock {\em IEEE Transactions on Information Theory}, 50(1):89--104, Jan
  2004.

\bibitem{LN2013}
B.~C. Levy and R.~Nikoukhah.
\newblock Robust state space filtering under incremental model perturbations
  subject to a relative entropy tolerance.
\newblock {\em IEEE Transactions on Automatic Control}, 58(3):682--695, March
  2013.

\bibitem{levy2013contraction}
B.~C. Levy and M.~Zorzi.
\newblock A contraction analysis of the convergence of risk-sensitive filters.
\newblock {\em SIAM J. Optimization Control}, 54(4):2154--2173, 2016.

\bibitem{iet5}
Shuai Liu, Yan Song, Guoliang Wei, and Xuegang Huang.
\newblock Rmpc-based security problem for polytopic uncertain system subject to
  deception attacks and persistent disturbances.
\newblock {\em IET Control Theory \& Applications}, 11(10):1611--1618, 2017.

\bibitem{iet4}
Yong-Hua Liu.
\newblock Saturated robust adaptive control for uncertain non-linear systems
  using a new approximate model.
\newblock {\em IET Control Theory \& Applications}, 11(6):870--876, 2017.

\bibitem{MACIEJOWSKI200922}
J.~Maciejowski.
\newblock Discussion on: âmin-max model predictive control of nonlinear
  systems: A unifying overview on stability".
\newblock {\em European Journal of Control}, 15(1):22--25, 2009.

\bibitem{bib:macie}
J.~M. Maciejowski.
\newblock {\em Predictive control: with constraints}.
\newblock Pearson education, 2001.

\bibitem{Mayne20142967}
D.~Q. Mayne.
\newblock Model predictive control: Recent developments and future promise.
\newblock {\em Automatica}, 50(12):2967 -- 2986, 2014.

\bibitem{bib:RMPCimp}
G.~De Nicolao, L.~Magni, and R.~Scattolini.
\newblock Robust predictive control of systems with uncertain impulse response.
\newblock {\em Automatica}, 32(10):1475--1479, 1996.

\bibitem{bib:RMPCH}
P.~E. Orukpe and I.~M. Jaimoukha.
\newblock Robust model predictive control based on mixed
  $\mathcal{H}_2/\mathcal{H}_\infty$ control approach.
\newblock {\em {Proceedings of the European Conference}}, pages 2223--2228,
  August 2009.

\bibitem{bib:surveyMPC}
S.~J. Qin and T.~A. Badgwell.
\newblock A survey of industrial model predictive control technology.
\newblock {\em {Control Engineering Practice II}}, pages 733--764, 2003.

\bibitem{iet1}
Ting Shi and Hongye Su.
\newblock Sampled-data mpc for lpv systems with input saturation.
\newblock {\em IET Control Theory \& Applications}, 8(17):1781--1788, 2014.

\bibitem{yang2015risk}
X.~Yang and J.~Maciejowski.
\newblock Risk-sensitive model predictive control with gaussian process models.
\newblock {\em IFAC-PapersOnLine}, 48(28):374--379, 2015.

\bibitem{YOON_2004}
M.~Yoon, V.~Ugrinovskii, and I.~Petersen.
\newblock Robust finite horizon minimax filtering for discrete-time stochastic
  uncertain systems.
\newblock {\em Syst. Control Lett.}, 52:99--112, 2004.

\bibitem{iet2}
Langwen Zhang, Jingcheng Wang, Yang Ge, and Bohui Wang.
\newblock Robust distributed model predictive control for uncertain networked
  control systems.
\newblock {\em IET Control Theory \& Applications}, 8(17):1843--1851, 2014.

\bibitem{BETA}
M.~Zorzi.
\newblock {A new family of high-resolution multivariate spectral estimators}.
\newblock {\em IEEE Trans. Autom. Control}, 59(4):892--904, Apr. 2014.

\bibitem{DUAL}
M.~Zorzi.
\newblock An interpretation of the dual problem of the {T}{H}{R}{E}{E}-like
  approaches.
\newblock {\em Automatica}, 62:87--92, 2015.

\bibitem{bib:RKdiv}
M.~Zorzi.
\newblock Multivariate spectral estimation based on the concept of optimal
  prediction.
\newblock {\em {IEEE Transactions on Automatic Control}}, 60(6):1647--1652,
  June 2015.

\bibitem{bib:RKweiner}
M.~Zorzi.
\newblock On the robustness of the {Bayes} and {W}iener estimators under model
  uncertainty.
\newblock {\em Automatica}, 83:133--140, 2017.

\bibitem{bib:rk}
M.~Zorzi.
\newblock Robust {K}alman filtering under model perturbations.
\newblock {\em IEEE Transactions on Automatic Control}, 62:2902--2907, Jun
  2017.

\bibitem{bib:RSconvergence}
M.~Zorzi and B.~C. Levy.
\newblock On the convergence of a risk sensitive like filter.
\newblock In {\em {54th IEEE Conference on Decision and Control}}, Osaka,Japan,
  December 2015.

\bibitem{convtau}
Mattia Zorzi.
\newblock Convergence analysis of a family of robust kalman filters based on
  the contraction principle.
\newblock {\em SIAM Journal on Optimization and Control}, 55(5):3116--3131,
  2017.

\end{thebibliography}
\end{document}